\documentclass{amsart}

\usepackage{amsmath, amsthm, amssymb, graphicx, tikz, zref} 
\usetikzlibrary{arrows}
\usepackage{enumitem}
\newtheorem{theorem}{Theorem}[section]

\newtheorem{conjecture}[theorem]{Conjecture}

\usepackage[margin=.6in]{geometry} 
\usepackage{tikz-cd}

\usepackage[backend=biber, 
style=alphabetic,
sorting=ynt
]{biblatex}
\usepackage{placeins}
\DeclareUnicodeCharacter{2212}{-}
\usepackage{comment}

\usepackage{pdfpages}
\usepackage{hyperref}
\usepackage{faktor}

\usepackage{mathabx,epsfig}

\hypersetup{
    colorlinks,
    linkcolor={red!50!black},
    citecolor={blue!50!black},
    urlcolor={blue!80!black}
}
\newcommand{\R}{\mathbb{R}}

\newcommand{\PP}{\mathbb{P}}\newcommand{\Z}{\mathbb{Z}}\newcommand{\C}{\mathbb{C}}

\usepackage{makecell}

\title{Exceptional Collections for Toric Fano Fourfolds}
\author{Jumari Querimit Ramirez, Yang (Hill) Zhang, Justin Son, and Reginald Anderson}
\date{May 2023}

\begin{document}

\begin{abstract}
Beilinson first gave a resolution of the diagonal for $\PP^n$. Generalizing this, a modification of the cellular resolution of the diagonal given by Bayer-Popescu- Sturmfels gives a (non-minimal, in general) virtual resolution of the diagonal for smooth projective toric varieties and toric Deligne-Mumford stacks which are a global quotient of a smooth projective variety by a finite abelian group. In the past year, Hanlon-Hicks-Lazarev gave in particular a symmetric, minimal resolution of the diagonal for smooth projective toric varieties. We give implications for exceptional collections on smooth projective toric Fano varieties in dimension 4. We find that for 72/124 smooth projective toric Fano 4-folds, the Hanlon-Hicks-Lazarev resolution of the diagonal yields a full strong exceptional collection of line bundles. Furthermore, we find that the success of the Hanlon-Hicks-Lazarev resolution to yield a full strong exceptional collection of line bundles coincides exactly with meeting a numerical criterion due to Bondal. \end{abstract}

\maketitle

\section{Introduction}

In this paper, we investigate exceptional collections for bounded derived category of coherent sheaves on smooth projective toric Fano fourfolds. Before explaining our results, let us explain our motivation for studying derived categories of toric varieties. First, King conjectured in unpublished notes that any smooth projective toric variety has a strong, full exceptional collection of line bundles. King's conjecture was proven false \cite{hille2006counterexample} \cite{Efimov2014}, and is currently checked on a case-by-case basis whether a given smooth projective toric variety has a strong, full exceptional collection of line bundles, and if so, how to construct such a collection. Though we work exclusively over $\C$, these questions have also been considered over more general fields \cite{ballard2020derivedcategoriescentrallysymmetricsmooth} \cite{BallardDuncanMcFaddin}. Additionally, full strong exceptional collections of line bundles on all 124 smooth projective toric Fano fourfolds in particular were constructed by Prabhu-Naik in \cite{PRABHUNAIK2017348} and used to create tilting bundles, building on work of Craw \cite{craw2008quiverrepresentationstoricgeometry} giving a new construction of toric varieties as a geometric quotient.

As a second motivation, note that through any two distinct points of $\PP^2(\C)$, the projective plane over the complex numbers, there exists a unique straight line, viewed as a rational degree $1$ curve as the image of $\PP^1(\C)$. Through five points in general position in $\PP^2(\C)$, there exists a unique conic. These kinds of enumerative problems have interested algebraic geometers since antiquity, by listing the number $n_d$ of degree $d$ rational curves through $3d-1$ points in general position in $\PP^2(\C)$. While $n_1=n_2=1$ have been known since antiquity, the invariants $n_3=12$ and $n_4=620$ say that the number of degree $d$ rational curves through $3d-1$ points in general position in $\PP^1$ is no longer $1$, which have been known since at least the 18th century \cite{zeuthen1873almindelige}. 
For the sake of exposition, we write $n_d$ in the formula below, though the enumerative invariants $n_d$ are actually related to the values $N_d$ in a recursive formula in a subtle way. The higher enumerative invariants $n_d$ for $\PP^2$ with $d\geq 5$ were discovered in the early 1990's, when they were given in Kontsevich-Manin's recursive formula \cite{Kontsevich_1994}

\[ n_d = \sum_{k+\ell = d} n_k n_\ell k^2 \ell \cdot \left[ \ell \binom{3d-4}{3k-2} - k \binom{3d-4}{3k-1} \right] \] 

using the theory of quantum cohomology and stable maps. Generalizing from $\PP^2(\C)$ to more abstract spaces, mathematicians and physicists alike have long been interested in enumerative invariants of spaces other than $\PP^2(\C)$. In 1879, Hermann Schubert found that a general quintic threefold contains 2875 lines \cite{Schubert1879}. In 1985, Katz found that there are 609250 conics on a general quintic threefold \cite{Katz1986}. After some back-and-forth between mathematicians and physicists \cite{Galison+2004+23+64}, it was determined that physicists Candelas-de la Ossa-Green-Parkes had correctly determined, $n_3$, the number of rational degree $3$ curves on the quintic threefold. To compute $n_3$, rather than studying the quintic threefold itself, Candelas-de la Ossa-Green-Parkes instead studied the mirror manifold $\check{Q}$ under the assumption that physics on either $Q$ or $\tilde{Q}$ gave similar physical theories. In his 1994 ICM Address, Kontsevich \cite{kontsevich1994homologicalalgebramirrorsymmetry} proposed that a categorical equivalence called Homological Mirror Symmetry (``HMS") should imply the enumerative mirror symmetry that Candelas-de la Ossa-Green-Parkes had used to calculate enumerative invariants on the quintic threefold. While originally stated for compact Calabi-Yau manifolds, HMS now since expanded to more general spaces, such as smooth projective toric Fano varieties. For a smooth projective toric Fano variety, one form of HMS for toric Fano varieties states that the bounded derived category of coherent sheaves on a smooth projective toric Fano variety $X$ should be equivalent to the Fukaya-Seidel category of the mirror Landau-Ginzburg model $\tilde{X} = ((\C^*)^n, W)$. This provides a second motivation to study the bounded derived category of coherent sheaves on a smooth projective toric Fano variety, as a presentation of $D^b_{Coh}(X)$ gives one side of the HMS correspondence in this setting. \\

A third motivation for studying derived categories comes from Dubrovin's conjecture 4.2.2 in \cite{dubrovin1998geometryanalytictheoryfrobenius}, which states that a Fano variety $X$ admits a full exceptional collection of length equal to $\text{dim} H^*(X)$ iff the quantum cohomology of $X$ is semisimple. Here, the full strong exceptional collections of line bundles for each smooth projective toric Fano fourfold given have length equal to $\text{dim }H^*(X)$, which is also the number of torus-fixed points in $X$ given by the number of vertices in the toric polytope $P$. \\

A best possible presentation of $D^b_{Coh}(X)$, the bounded derived category of coherent sheaves on a smooth projective toric variety $X$, comes from the presence of a full strong exceptional collection of objects $\mathcal{E}$ \cite{BondalRepnAssocAlg}. For a smooth projective Fano or general type variety, $D^b_{Coh}(X)$ determines $X$ up to isomorphism \cite{bondal2002derivedcategoriescoherentsheaves}. Beilinson's resolution of the diagonal \cite{Beilinson1978} gave a full strong exceptional collection of line bundles for $\PP^n(\C)$, which is a smooth projective toric Fano variety. Generalizing Beilinson's resolution of the diagonal, Bayer-Popescu-Sturmfels \cite{bayer-popescu-sturmfels} gave a cellular resolution of the diagonal for a proper subclass of smooth projective toric varieties which they called ``unimodular," which is a more restrictive condition than being nonsingular. Generalizing this, a non-minimal virtual cellular resolution of the diagonal was given by the fourth author in \cite{anderson2023resolutiondiagonalsmoothtoric} for smooth projective toric varieties and in \cite{anderson2023resolutiondiagonaltoricdelignemumford} for toric Deligne-Mumford stacks which are a global quotient of a smooth projective toric variety by a finite abelian group. In \cite{hanlon2023resolutions}, Hanlon-Hicks-Lazarev gave a minimal virtual resolution of the diagonal for any toric subvariety of a smooth projective toric variety. Here, we investigate the question: ``For which smooth projective toric Fano fourfolds does the Hanlon-Hicks-Lazarev resolution of the diagonal yield a full strong exceptional collection of line bundles?" Our paper focuses on the following Theorem~\ref{thm: Thm1}.

\begin{theorem} \label{thm: Thm1} The Hanlon-Hicks-Lazarev resolution of the diagonal yields a full strong exceptional collection of line bundles for 72 out of 124 smooth projective toric Fano fourfolds. Furthermore, the success of the Hanlon-Hicks-Lazarev resolution of the diagonal to yield a full strong exceptional collection of line bundles coincides exactly with Bondal's numerical criterion being successfully met. \end{theorem}

Smooth projective toric varieties for which the Hanlon-Hicks-Lazarev resolution yields a full strong exceptional collection of line bundles are referred to as Bondal-Ruan type \cite{favero2022homotopy} \cite{bondalOberwolfach}.

\newpage 

\begin{center}
\includegraphics[width=.9\textwidth]{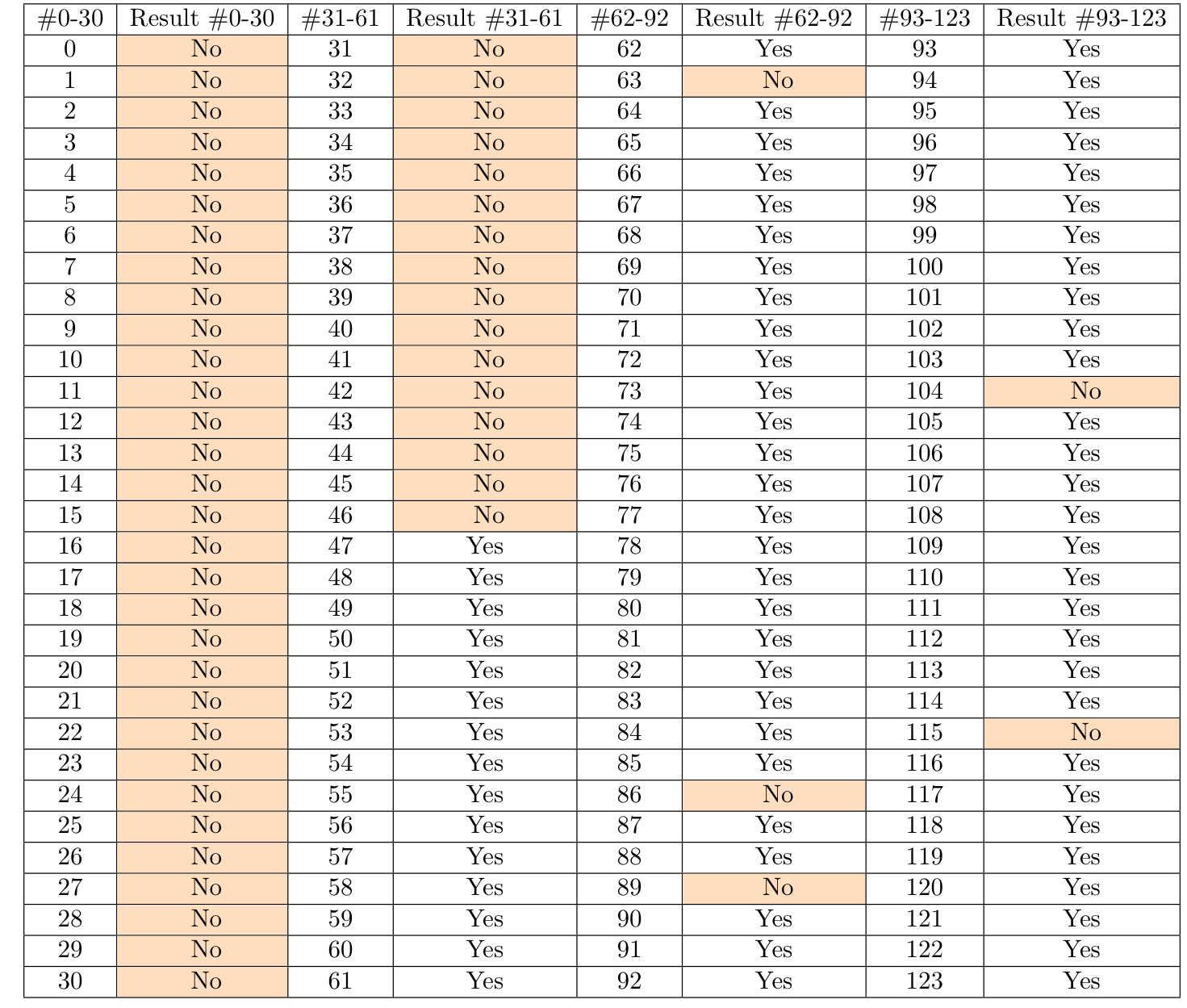}
\end{center}

\section{Acknowledgements}

We are grateful to Jay Yang for providing Macaulay2 code reproducing the Hanlon-Hicks-Lazarev resolution, with contributions to the code from Mahrud Sayrafi. We are also grateful for discussions with Daniel Erman, Christine Berkesch, Jesse Huang, and Andrew Hanlon. The second and fourth author were also partially supported from the Claremont McKenna College Summer Research Program for this project.

\section{Background}

We work over the algebraically closed field $k=\overline{k}=\C$. Let $\mathcal{C}$ denote a k-linear triangulated category. Recall \cite{huybrechts2006fourier} that in a $k$-linear triangulated category $\mathcal{D}$, an object $E$ is called \textbf{exceptional} if \[ \text{Hom}_\mathcal{D}(E, E[\ell]) = \begin{cases} k & \ell = 0, \\ 0 & \text{ else.} \end{cases} \] 

An \textbf{exceptional sequence} $\mathcal{E}$ is a sequence $E_1, \dots, E_r$ of exceptional objects such that $Hom_\mathcal{E}(E_i, E_j[\ell]) = 0$ for $i>j$ and all $\ell$. That is, 

\[ \text{Hom}_\mathcal{D}(E_i, E_j[\ell]) = \begin{cases} k & i=j, \ell=0 \\ 0 & \text{ if } i>j \text{ or if } \ell\neq 0, i=j. \end{cases} \] 

An exceptional sequence is \textbf{full} if $\mathcal{D}$ is generated by $\{E_i\}$. An exceptional collection $E_1, \dots, E_r$ is \textbf{strong} if \[ \text{Hom}_{\mathcal{D}}(E_i, E_j[\ell]) = 0\] for all $i,j$, and $\ell\neq 0$. Here, $\mathcal{D} = D^b_{Coh}(X)$, the bounded derived category of coherent sheaves on a smooth projective toric Fano variety of dimension $4$ over $\mathbb{C}$. Given a full strong exceptional collection of line bundles $\mathcal{E} = \{E_i\}_{i=1}^r$ on a smooth projective variety, we form $\mathcal{G} = \bigoplus_{i=1}^r E_i$ and $A =\text{End}_{D^b_{Coh}(X)}(\mathcal{G})$ so that $D^b_{Coh}(X) \simeq D^b(\text{mod}-A)$ \cite{BondalRepnAssocAlg}. 

\section{Methods}

The methods that we have used for this project rely on the polymake database of 4-dimensional Fano polytopes \cite{Assarf_Gawrilow_Herr_Joswig_Lorenz_Paffenholz_Rehn_2016} for the corresponding classification of smooth projective toric Fano fourfolds. The polymake database uses an algorithm due to \O bro for the classification of smooth reflexive 4-polytopes up to isomorphism \cite{øbro2007algorithmclassificationsmoothfano}. Here, the fan $\Sigma$ for the toric variety $X$ has rays given by inner-pointing normals to the polytope $P$ from the polymake database \cite{C-L-S}. The convention from the polymake database for the half-space representation of the polytope $P$ is that each row $(a_0, \dots a_d)$ encodes the linear inequality $a_0 + a_1x_1 + \cdots + a_d x_d \geq 0$. Macaulay2 uses a different half-space representation convention in the command \begin{verbatim} polyhedronFromHData(I, v) \end{verbatim} which corresponds to $\{x\in \R^d \text{ }|\text{ } I \cdot x \leq v\}$ where $I$ is an $N$-by-d matrix for $N$ the number of defining inequalities of the polytope $P$, and $v$ a column vector in $\R^N$, so we necessarily translate between the two conventions. Given the polymake polytope $P$, we apply the Hanlon-Hicks-Lazarev resolution to $\mathcal{O}_\Delta$ for $X \cong \Delta \subset X\times X$ the diagonal embedding in terms of line bundles. Since the Hanlon-Hicks-Lazarev resolution is symmetric in terms of which line bundles appear on either of the left- or right-hand side of $X\times X$, we consider the collection of line bundles $\mathcal{E}$ appearing on the left-hand side and ask whether there exists an ordering of the collection of line bundles $\mathcal{E}$ for which $\mathcal{E}$ give a full exceptional collection, and further, whether this collection is strong. Here, there is a dichotomy in the fact that if there exists an ordering of the objects $\{E_i\}$ in $\mathcal{E}$, then this collection is also strong; and if ever a non-zero $Hom(E_i, E_j[\ell])$ appears with $\ell\neq 0$, then there does not exist an ordering for which $\mathcal{E}$ gives an exceptional collection. We do not speak of the collection $\mathcal{E}$ being strong if $\mathcal{E}$ is not also exceptional. 

Here, the collection $\mathcal{E}$ will always be full in virtue of giving a locally-free resolution of $\mathcal{O}_\Delta$ in terms of line bundles. The obstruction to the existence of an ordering for which $\mathcal{E}$ is exceptional is exactly the presence of a directed cycle of positive length in the directed quiver which lists $ \text{Hom}_{D^b_{Coh}(X)}(E_i, E_j[\ell])$ 
for all $i,j$, and $\ell$. To compute $Hom_{D^b_{Coh}(X)}(E_i, E_j[\ell])$, we make frequent use of the fact that for line bundles $\mathcal{O}(D)$ and $\mathcal{O}(E)$ on a smooth projective toric variety, 

\[ \text{Hom}_{D^b_{Coh}(X)}(\mathcal{O}(D), \mathcal{O}(E)) \cong H^*(X, \mathcal{O}(E-D) ). \]

When there exists an ordering for which $\mathcal{E}$ is exceptional (which here also implies that $\mathcal{E}$ is strong), we indicate this by giving a lower-triangular matrix whose $(i,j)$ entry indicates the rank of $Hom_{D^b_{Coh}(X)}(E_j, E_i[0])$ and a directed quiver as in \cite{BondalRepnAssocAlg}. 

The fourth author constructed a virtual cellular resolution of the diagonal for smooth projective toric varieties in \cite{anderson2023resolutiondiagonalsmoothtoric} 
generalizing Bayer-Popescu-Sturmfels' approach \cite{bayer-popescu-sturmfels}, which in turn generalized Beilinson's resolution of the diagonal for $\PP^n$ \cite{Beilinson1978}. Smooth projective toric varieties have a presentation of the class group over integer-valued matrices via the fundamental exact sequence \[ 0 \rightarrow M \stackrel{B}{\rightarrow} \Z^{|\mathcal{A}|} \stackrel{\pi}{\rightarrow} \text{Cl}(X) \rightarrow 0 \] where $M$ is the character lattice of the Zariski-dense open algebraic torus $T \cong (\C^*)^m \subseteq X$ of rank $m$, $\mathcal{A}=\Sigma(1)$ is the set of primitive ray generators with $|\Sigma(1)|=n$, and, fixing a basis for $M$, $B$ is represented by the $n \times m$ matrix with rows as primitive ray generators, which are the elements of $\mathcal{A}$. Let $L$ denote the image of $B$ in $\Z^{|\mathcal{A}|}$, and consider $\R L = L\otimes_{\Z} \R \subset \R^{|\mathcal{A}|} = \Z^{|\mathcal{A}|} \otimes_\Z \R$. $\R L$ and the infinite hyperplane arrangement $\mathcal{H} = \{ x_i = j \text{ }|\text{ } 1 \leq i \leq n, j \in \Z\} \subset \R^{|\mathcal{A}|}$ were used in \cite{anderson2023resolutiondiagonalsmoothtoric} to construct an infinite and quotient cellular complex, respectively, to build a cellular resolution of the diagonal. Hanlon-Hicks-Lazarev \cite{hanlon2023resolutions} use the same hyperplane arrangement and a different convention on Laurent-monomial labelings for the finite quotient cellular complex which agrees with \cite{bondalOberwolfach} to build a minimal resolution of the diagonal\cite{brown2023short} for smooth projective toric varieties. \\

In \cite{bondalOberwolfach}, Bondal gave a criterion for a certain map $\Phi$ which is used in the construction of the Hanlon-Hicks-Lazarev resolution of the diagonal to yield a full strong exceptional collection of line bundles on the smooth projective toric variety $X$ of dimension $n$: for each irreducible toric curve $C$, denote by $(D_1, \dots, D_{n-1})$ the irreducible toric divisors that contain $C$ and by $(a_1, \dots, a_{n-1})$ the corresponding intersection numbers with $C$. Bondal's criterion is that all $a_i \geq -1$, with $a_i=-1$ appearing no more than once. We refer to this condition, which implies that the Hanlon-Hicks-Lazarev resolution of the diagonal yields a full strong exceptional collection of line bundles when it is satisfied, simply as Bondal's numerical criterion. Through direct calculation, we also verify that failing Bondal's numerical criterion coincides exactly with the success of the Hanlon-Hicks-Lazarev resolution of the diagonal to yield a full, strong, exceptional collection of line bundles on a smoth projective toric variety in dimension 4. 

\section{Data collection}
Our data for checking whether the Hanlon-Hicks-Lazarev resolution of the diagonal yields a full strong exceptional collection on each of the 124 smooth projective toric Fano 4-folds is available at \url{https://github.com/reggiea91/CMC-Fourfolds.git}. For the sake of exposition, we include a negative and positive example, respectively, from our data collection procedure. 

\subsection{Polytope F.4D.0000}

Let $P$ denote the polytope F.4D.0000 in polymake with half-space representation:
\[ \left[ \begin{matrix} 1 & 0& -1 &1 &0 \\
1& 0& -1& 1& 0 \\
1 &0 &0& 0& -1 \\
1 &-1 &0 &0& 0 \\
1& 0& 0 &1 &0 \\
1 &1 &1 &-3 &1 \end{matrix} \right] \]
with vertices given by the columns of

\[  \left( \begin{matrix} 1&1&1&-3&1&-5&1&1&1&0&1 \\1&0&1&1&1&0&0&-6&0&1&1 \\ 1&-1&0&0&0&-1&-1&-1&1&1&1 \\0&-5&-3&1&1&1&1&1&1&1&1 \end{matrix} \right). \] Let X denote the complete toric variety associated to P, with rays

\[ \{ (-1,0,0,0), (0,-1,0,0), (0,0,-1,0), (0,0,1,0), (0,-1,1,0), (0,0,0,-1), (1,1,-3,1) \} \]
and presentation of the class group

\[  \Z^7 \stackrel{ \left( \begin{matrix} 0&0&1&1&0&0&0 \\ 0&-1&1&0&1&0&0 \\ 1&1&-3&0&0&1&1 \end{matrix} \right) }{\longrightarrow} \Z^3. \]
The Hanlon-Hicks-Lazarev resolution yields the free ranks (written as an ungraded resolution of S-modules, for S the
homogeneous coordinate ring of $Y \cong X \times X)$: 

\[0 \rightarrow S^{10} \rightarrow S^{36} \rightarrow S^{47} \rightarrow S^{25} \rightarrow S^4 \rightarrow 0. \]

The Hanlon-Hicks-Lazarev resolution yields the collection of line bundles $\{\mathcal{O}(a_1, a_2)\}$ for $(a_1, a_2)$ appearing in

\begin{align*} 
\mathcal{E}  = &\{  (-1,-1,1), (0,0,0),(-1,-1,2), (-1,0,1), (-1,0,0), (-1,-1,0), (0,0,-1),(-1,-1,-1), \dots \\
& \dots (-1,0,-1), (0,0,-2), (0,0,-3),(-1,-1,-2), (-1,0,-3) \}  \end{align*} 
which appear on one side of the Hanlon-Hicks-Lazarev resolution for which the collection is exceptional, since there are a pair of line bundles $\mathcal{O}(-1, 0, -3)$ and $\mathcal{O}(0,-1,4)$ with non-zero Hom's in both directions. That is,

\[Hom^\bullet_{D^b(X)}( \mathcal{O}(-1, 0, -3), \mathcal{O}(0,-1,4) ) =  \begin{cases} \C & \text{ in degree }3\\
0 & \text{ else}. \end{cases}  \] 

and

\[ Hom^\bullet_{D^b(X)}( \mathcal{O}(0,-1,4), \mathcal{O}(-1, 0, -3)) =  \begin{cases} \C^{20} &\text{ in degree }0 \\ 0 & \text{ else}. \end{cases} \]

Here, we see that there does \textbf{not} exist an ordering of $\mathcal{E}$ for which we get an exceptional collection of line bundles. \\


Furthermore, we verify that Bondal's numerical criterion fails for this variety, since the torus-invariant curve corresponding to the 3-cone with primtitive ray generators $\{ (-1,  0, 0, 0), ( 0, -1, 0, 0),( 0, -1, 1, 0) \}$ intersects the torus-invariant divisor corresponding to primitive ray generator $(0, -1, 0, 0)$ in $-2$.

\subsection{Polytope F.4D.0085 } Let $P$ denote the polytope F.4D.0085 in the polymake database, with half-space representation

\[ \left[ \begin{matrix} 1 & 0 & 0 & 0 & -1 \\
1 & 0 & -1 & 0 & 0 \\
1 & -1 & 0 & 0 & 0 \\
1 & 0 & -1 & 1 & 0 \\
1 & 0 & 0 & -1 & 0 \\
1 & 0 & 1 & 0 & 1 \\
1 & 1 & 1 & 0 & 1 \end{matrix} \right] \] 

and vertices of $P$ given by the columns of:
\[ \left( \begin{matrix}0 & 1 & 0 & 1 & -3 & 1 & 0 & 1 & 0 & 1 & -3 & 1 \\
1 & 1 & 1 & 1 & 1 & 1 & -2 & -2 & -2 & -2 & 1 & 1 \\
0 & 0 & 0 & 1 & 0 & 0 & -3 & -3 & 1 & 1 & 1 & 1 \\
-2 & -2 & -2 & -2 & 1 & 1 & 1 & 1 & 1 & 1 & 1 & 1 \end{matrix} \right) \]

Let $X$ denote the complete toric variety associated to $P$, with rays

\[ \{(-1, 0, 0, 0), (0, -1, 0, 0), (0, 0, -1, 0), (0, 1, 1, 0), (0, 0, 0, -1), (0, 1, 0, 1), (1, 1, 0, 1)\} \]

We use as presentation of $Cl(X)$:

\[ \Z^7 \stackrel{\left( 
\begin{matrix}
0 & -1 & 1 & 1 & 0 & 0 & 0 \\
0 & 1 & 0 & 0 & 1 & 1 & 0 \\
1 & 1 & 0 & 0 & 1 & 0 & 1  \end{matrix} \right) }{\longrightarrow } \Z^3. \] 

The Hanlon-Hicks-Lazarev resolution yields the free ranks (written as an ungraded resolution of S-modules, for $S$ the homogeneous coordinate ring of $Y \cong X \times X)$:

\[ 0 \rightarrow S^8 \rightarrow S^{24} \rightarrow S^{25} \rightarrow S^{10} \rightarrow S^1 \rightarrow 0.\]

There are 12 line bundles appearing on left-hand side which give a full strong exceptional collection of line bundles as follows:
\begin{align*} \mathcal{E} =\{ & (-1, -2, -3), (-1, -2, -2), (0, -2, -3), (-1, -1, -2), (0, -2, -2), (0, -1, -2), \dots \\
& \dots(-1, -1, -1), (-1, 0, -1), (0, 0, -1), (-1, 0, 0), (0, -1, -1), (0, 0, 0)
 \} \end{align*}

Here $Hom_{D^b_{Coh}(X)}(E_i, E_j)$ is concentrated in degree 0 for all $i$ and $j$, with the rank of $Hom^0(E_j, E_i)$ given by the $(i,j)$ entry of the following matrix: 

\begin{align*}
\left( \begin{matrix} 
1 & 0 & 0 & 0 & 0 & 0 & 0 & 0 & 0 & 0 & 0 & 0 \\
2 & 1 & 0 & 0 & 0 & 0 & 0 & 0 & 0 & 0 & 0 & 0 \\
2 & 0 & 1 & 0 & 0 & 0 & 0 & 0 & 0 & 0 & 0 & 0 \\
5 & 1 & 1 & 1 & 0 & 0 & 0 & 0 & 0 & 0 & 0 & 0 \\
4 & 2 & 2 & 0 & 1 & 0 & 0 & 0 & 0 & 0 & 0 & 0 \\
9 & 2 & 5 & 2 & 1 & 1 & 0 & 0 & 0 & 0 & 0 & 0 \\
9 & 5 & 2 & 2 & 1 & 0 & 1 & 0 & 0 & 0 & 0 & 0 \\
15 & 5 & 5 & 5 & 1 & 1 & 1 & 1 & 0 & 0 & 0 & 0 \\
16 & 9 & 9 & 4 & 5 & 2 & 2 & 0 & 1 & 0 & 0 & 0 \\
25 & 15 & 9 & 9 & 5 & 2 & 5 & 2 & 1 & 1 & 0 & 0 \\
25 & 9 & 9 & 9 & 5 & 5 & 2 & 2 & 1 & 0 & 1 & 0 \\
41 & 25 & 25 & 16 & 15 & 9 & 9 & 4 & 5 & 2 & 2 & 1
\end{matrix} \right) 
\end{align*}

The directed quiver $\mathcal{Q}$ showing nonzero $Hom^0(E_j, E_i)$, for $i\neq j$ is then given by: 

\FloatBarrier
\begin{figure}[h]
\includegraphics{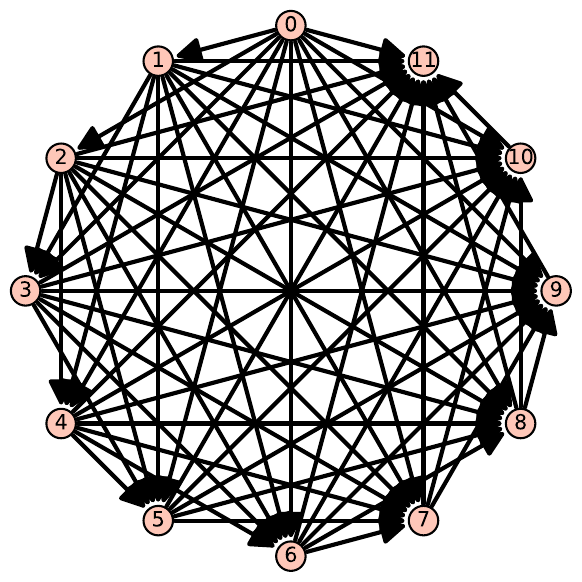}
\label{fig: F.4D.0085}
\caption{Quiver showing nonzero $Hom^0(E_i, E_j)$ for $i\neq j$ in $\mathcal{E}$. Indexing on vertices from $0$ to $11=|\mathcal{E}|-1$.} 
\end{figure}
\FloatBarrier

\section{Future directions} 
In dimension 1, the Hanlon-Hicks-Lazarev resolution yields a full strong exceptional collection of line bundles for the unique smooth projective toric Fano curve. In this case, this coincides with Beilinson's collection for $\PP^1$. In dimension 2, the resolution of the diagonal developed by the fourth author coincides with the Hanlon-Hicks-Lazarev resolution of the diagonal for all 5 smooth projective toric Fano surfaces, since all are unimodular in the sense of Bayer-Popescu-Sturmfels. In this case, either resolution of the diagonal yields a full strong exceptional colleciton for all 5 smooth projective toric Fano surfaces \cite{anderson2024exceptionalcollectionslinebundles}. In dimension 3, 2 out of 18 smooth projective toric Fano 3-folds provably do not admit a full strong exceptional collection of line bundles from the Hanlon-Hicks-Lazarev resolution of the diagonal due to failing Bondal's numerical condition \cite{bondalOberwolfach}. We have found that the same is true for all 52 out of 124 smooth projective toric Fano fourfolds for which the Hanlon-Hicks-Lazarev resolution of the diagonal fails to yield a full strong exceptional collection of line bundles. This leads to the following conjecture.

\begin{conjecture}

Let $f(d): \Z_{>0} \rightarrow [0,1]$ denote the proportion of smooth projective toric Fano d-folds for which the Hanlon-Hicks-Lazarev resolution of the diagonal yields a full strong exceptional collection of line bundles. Then we have 

\center 
   
  \[ \begin{array}{|c|c|} 
   \hline 
    d & f(d)\\
    \hline &  \\
    1 & 1 \\
    \hline &  \\
    2 & 1 = \frac{5}{5} \\
    \hline & \\
    3 & \frac{16}{18} \\
\hline & \\
    4 & \frac{72}{124}\\ 
    \hline  \end{array} \]

\flushleft 
We conjecture that $f(d)$ is strictly monotone decreasing for $d>1$, and that $\lim_{d\rightarrow \infty} f(d) \neq 0$. Furthermore, we conjecture that the success of the Hanlon-Hicks-Lazarev resolution of the diagonal to yield a full strong exceptional collection of line bundles coincides exactly with passing Bondal's numerical criterion in any dimension $d\geq 5$. 

\end{conjecture}

\clearpage
\printbibliography

@misc{anderson2024exceptionalcollectionslinebundles,
      title={Exceptional Collections of Line Bundles for Smooth Toric Fano Surfaces and Threefolds}, 
      author={Reginald Anderson},
      year={2024},
      eprint={2403.09663},
      archivePrefix={arXiv},
      primaryClass={math.AG},
      url={https://arxiv.org/abs/2403.09663}, 
}

@misc{craw2008quiverrepresentationstoricgeometry,
      title={Quiver representations in toric geometry}, 
      author={Alastair Craw},
      year={2008},
      eprint={0807.2191},
      archivePrefix={arXiv},
      primaryClass={math.AG},
      url={https://arxiv.org/abs/0807.2191}, 
}

@article{PRABHUNAIK2017348,
title = {Tilting bundles on toric Fano fourfolds},
journal = {Journal of Algebra},
volume = {471},
pages = {348-398},
year = {2017},
issn = {0021-8693},
doi = {https://doi.org/10.1016/j.jalgebra.2016.09.007},
url = {https://www.sciencedirect.com/science/article/pii/S0021869316303118},
author = {Nathan Prabhu-Naik}
}

@article{BallardDuncanMcFaddin,
author = {Matthew Ballard and Alexander Duncan and Patrick McFaddin},
title = {{On derived categories of arithmetic toric varieties}},
volume = {4},
journal = {Annals of K-Theory},
number = {2},
publisher = {MSP},
pages = {211 -- 242},
year = {2019},
doi = {10.2140/akt.2019.4.211},
URL = {https://doi.org/10.2140/akt.2019.4.211}
}

@misc{ballard2020derivedcategoriescentrallysymmetricsmooth,
      title={Derived categories of centrally-symmetric smooth toric Fano varieties}, 
      author={Matthew R Ballard and Alexander Duncan and Patrick K. McFaddin},
      year={2020},
      eprint={1812.09392},
      archivePrefix={arXiv},
      primaryClass={math.AG},
      url={https://arxiv.org/abs/1812.09392}, 
}

@misc{øbro2007algorithmclassificationsmoothfano,
      title={An algorithm for the classification of smooth Fano polytopes}, 
      author={Mikkel Øbro},
      year={2007},
      eprint={0704.0049},
      archivePrefix={arXiv},
      primaryClass={math.CO},
      url={https://arxiv.org/abs/0704.0049}, 
}

@misc{bondal2002derivedcategoriescoherentsheaves,
      title={Derived categories of coherent sheaves}, 
      author={Alexei Bondal and Dmitri Orlov},
      year={2002},
      eprint={math/0206295},
      archivePrefix={arXiv},
      primaryClass={math.AG},
      url={https://arxiv.org/abs/math/0206295}, 
}

@misc{anderson2023resolutiondiagonaltoricdelignemumford,
      title={A Resolution of the Diagonal for Toric Deligne-Mumford Stacks}, 
      author={Reginald Anderson},
      year={2023},
      eprint={2303.17497},
      archivePrefix={arXiv},
      primaryClass={math.AG},
      url={https://arxiv.org/abs/2303.17497}, 
}

@misc{anderson2023resolutiondiagonalsmoothtoric,
      title={A Resolution of the Diagonal for Smooth Toric Varieties}, 
      author={Reginald Anderson},
      year={2023},
      eprint={2403.09653},
      archivePrefix={arXiv},
      primaryClass={math.AG},
      url={https://arxiv.org/abs/2403.09653}, 
}

@misc{dubrovin1998geometryanalytictheoryfrobenius,
      title={Geometry and analytic theory of Frobenius manifolds}, 
      author={Boris Dubrovin},
      year={1998},
      eprint={math/9807034},
      archivePrefix={arXiv},
      primaryClass={math.AG},
      url={https://arxiv.org/abs/math/9807034}, 
}

@misc{kontsevich1994homologicalalgebramirrorsymmetry,
      title={Homological Algebra of Mirror Symmetry}, 
      author={Maxim Kontsevich},
      year={1994},
      eprint={alg-geom/9411018},
      archivePrefix={arXiv},
      primaryClass={alg-geom},
      url={https://arxiv.org/abs/alg-geom/9411018}, 
}

@inbook{Galison+2004+23+64,
url = {https://doi.org/10.1515/9780822390084-002},
title = {1 Mirror symmetry: persons, values, and objects},
booktitle = {Growing Explanations},
booktitle = {Historical Perspectives on Recent Science},
author = {Peter Galison},
editor = {M. Norton Wise and Barbara Herrnstein Smith and E. Roy Weintraub},
publisher = {Duke University Press},
address = {New York, USA},
pages = {23--64},
doi = {doi:10.1515/9780822390084-002},
isbn = {9780822390084},
year = {2004},
lastchecked = {2024-07-19}
}

@article{Katz1986,
author = {Katz, Sheldon},
journal = {Compositio Mathematica},
language = {eng},
number = {2},
pages = {151-162},
publisher = {Martinus Nijhoff Publishers},
title = {On the finiteness of rational curves on quintic threefolds},
url = {http://eudml.org/doc/89802},
volume = {60},
year = {1986},
}

@book{Schubert1879,
  author = {Hermann Schubert},
  year = {1879},
  title = {Kalkül der abzählenden Geometrie},
  publisher = {Leipzig B.G. Teubner}
}

@article{Kontsevich_1994,
   title={Gromov-Witten classes, quantum cohomology, and enumerative geometry},
   volume={164},
   ISSN={1432-0916},
   url={http://dx.doi.org/10.1007/BF02101490},
   DOI={10.1007/bf02101490},
   number={3},
   journal={Communications in Mathematical Physics},
   publisher={Springer Science and Business Media LLC},
   author={Kontsevich, M. and Manin, Yu.},
   year={1994},
   month=aug, pages={525–562} }

@book{zeuthen1873almindelige,
  title={Almindelige Egenskaber ved Systemer af plane Kurver, med Anvendelse til Bestemmelse af Karakteristikerne i de element{\ae}re Systemer af fjerde Orden},
  author={Zeuthen, H.G.},
  number={v. 1},
  series={Almindelige Egenskaber ved Systemer af plane Kurver, med Anvendelse til Bestemmelse af Karakteristikerne i de element{\ae}re Systemer af fjerde Orden},
  url={https://books.google.com/books?id=n4taQwAACAAJ},
  year={1873},
  publisher={Luno}
}

@article{BondalRepnAssocAlg,
doi = {10.1070/IM1990v034n01ABEH000583},
url = {https://dx.doi.org/10.1070/IM1990v034n01ABEH000583},
year = {1990},
month = {feb},
publisher = {},
volume = {34},
number = {1},
pages = {23},
author = {A I Bondal},
title = {REPRESENTATION OF ASSOCIATIVE ALGEBRAS
AND COHERENT SHEAVES},
journal = {Mathematics of the USSR-Izvestiya}
}

@article{bondalOberwolfach,
author = {Bondal, Alexey},
year = {2006},
month = {01},
pages = {284-286},
title = {Derived categories of toric varieties},
volume = {3},
journal = {Oberwolfach Reports}
}

@article{Assarf_Gawrilow_Herr_Joswig_Lorenz_Paffenholz_Rehn_2016, title={Computing convex hulls and counting integer points with polymake}, volume={9}, DOI={10.1007/s12532-016-0104-z}, number={1}, journal={Mathematical Programming Computation}, author={Assarf, Benjamin and Gawrilow, Ewgenij and Herr, Katrin and Joswig, Michael and Lorenz, Benjamin and Paffenholz, Andreas and Rehn, Thomas}, year={2016}, month={May}, pages={1–38}}

@misc{brown2023short,
      title={A short proof of the Hanlon-Hicks-Lazarev Theorem}, 
      author={Michael K. Brown and Daniel Erman},
      year={2023},
      eprint={2303.14319},
      archivePrefix={arXiv},
      primaryClass={math.AG}
}

@misc{favero2022homotopy,
      title={Homotopy Path Algebras}, 
      author={David Favero and Jesse Huang},
      year={2022},
      eprint={2205.03730},
      archivePrefix={arXiv},
      primaryClass={math.AG}
}

@misc{hanlon2023resolutions,
      title={Resolutions of toric subvarieties by line bundles and applications}, 
      author={Andrew Hanlon and Jeff Hicks and Oleg Lazarev},
      year={2023},
      eprint={2303.03763},
      archivePrefix={arXiv},
      primaryClass={math.AG}
}

@misc{hille2006counterexample,
      title={A Counterexample to King's Conjecture}, 
      author={Lutz Hille and Markus Perling},
      year={2006},
      eprint={math/0602258},
      archivePrefix={arXiv},
      primaryClass={math.AG}
}

@article{Efimov2014,
	doi = {10.1112/jlms/jdu037},
	url = {https://doi.org/10.1112%2Fjlms%2Fjdu037},
	year = 2014,
	month = {jul},
	publisher = {Wiley},
	volume = {90},
	number = {2},
	pages = {350--372},
	author = {Alexander I. Efimov},
	title = {Maximal lengths of exceptional collections of line bundles},
	journal = {Journal of the London Mathematical Society}
}

@book{C-L-S,
title = "Toric varieties",
author = "Cox, {David A.} and Little, {John B.} and Schenck, {Henry K.}",
year = "2011",
doi = "10.1090/gsm/124",
language = "English (US)",
isbn = "978-0-8218-4819-7",
volume = "124",
series = "Graduate Studies in Mathematics",
publisher = "American Mathematical Society",
address = "United States",
}

@misc{bayer-popescu-sturmfels,
  doi = {10.48550/ARXIV.MATH/9912247},
  
  url = {https://arxiv.org/abs/math/9912247},
  
  author = {Bayer, Dave and Popescu, Sorin and Sturmfels, Bernd},
  
  title = {Syzygies of Unimodular Lawrence Ideals},
  
  publisher = {arXiv},
  
  year = {1999},
  
  copyright = {Assumed arXiv.org perpetual, non-exclusive license to distribute this article for submissions made before January 2004}
}

@article{Beilinson1978,

author = {Beilinson, A. A.},
title = {Coherent Sheaves on Pn and Problems of Linear Algebra},
volume = {12},
journal = {Functional Analysis and Its Applications},
number = {3},
pages = {214 -- 216},
year = {1978},
doi = {10.1007/BF01681436},
url = {https://doi.org/10.1007/BF01681436}

}

@book{huybrechts2006fourier,
  title={Fourier-Mukai Transforms in Algebraic Geometry},
  author={Huybrechts, D.},
  isbn={9780199296866},
  lccn={2006298244},
  series={Oxford Mathematical Monographs},
  url={https://books.google.com/books?id=9HQTDAAAQBAJ},
  year={2006},
  publisher={Clarendon Press}
}

\end{document}